\documentstyle[12pt,amssymb]{amsart}

\ifx\mathrm\undefined\let\mathrm\bf\fi
\ifx\mathbf\undefined\let\mathbf\bf\fi
\ifx\mathfrak\undefined\let\mathfrak\frak\fi
\ifx\mathcal\undefined\let\mathcal\cal\fi
\ifx\mathbb\undefined\let\mathbb\Bbb\fi
\ifx\emph\undefined\let\emph\it\fi
\newcommand{\Z}{{\mathbb Z}}

\newcommand{\R}{{\mathbb R}}
\newcommand{\C}{{\mathbb C}}
\newcommand{\Q}{{\mathbb Q}}
\newcommand{\Pee}{{\mathbb P}}
\newcommand{\Ref}[1]{{$($\ref{#1}$)$}}
\newcommand{\be}{\begin{displaymath}}
\newcommand{\ee}{\end{displaymath}}
\newcommand{\bea}{\begin{eqnarray*}}
\newcommand{\eea}{\end{eqnarray*}}
\newcommand{\tr}{{\mathrm{tr}}}
\newcommand{\T}{{\mathcal{T}}}
\newcommand{\Mu}{{\mathrm{M}}}
\newenvironment{prf}{\noindent{\it Proof\/}:}{$\;\Box$
\par\medskip}

\newtheorem%
{thm}{Theorem}[section]
\newtheorem%
{proposition}[thm]{Proposition}
\newtheorem%
{lemma}[thm]{Lemma}
\newtheorem%
{lemmadef}[thm]{Lemma-Definition}
\newtheorem%
{corollary}[thm]{Corollary}
\newtheorem%
{conjecture}[thm]{Conjecture}
\newenvironment{definition}
{\noindent{\bf Definition\/}:}{\par\medskip}

\title{Algebraic integrability of the two-body Ruijsenaars operator}
\author[G. Felder and A. Varchenko]{Giovanni Felder${}^{*,1}$ 
\and Alexander Varchenko${}^{**,2}$}
\thanks{${}^1$Permanent address: D-MATH, ETH-Zentrum, 8092
Z\"urich, Switzerland}
\thanks{${}^2$Supported in part by NSF grant  DMS-9501290}
\begin{document}
\maketitle
\medskip
\centerline{\it ${}^*$Centre Emile Borel, Institut Henri Poincar\'e,
Universit\'e Pierre et Marie Curie,}
\centerline{\it 75231 Paris Cedex 05, France}
\medskip
\centerline{\it ${}^{**}$Department of Mathematics,
University of North Carolina at Chapel Hill,}
\centerline{\it Chapel Hill, NC 27599-3250, USA}
\medskip
\centerline{October 1996}
\begin{abstract} We study the algebra of difference operators that
commute with the two-body Ruijsenaars operator, a $q$-deformation
of the Lam\'e differential operator, for generic values of the deformation
parameter. The algebra is commutative. It
is the algebra of polynomial functions on an affine hyperelliptic
curve $Y^2=P(X^2)$. We also compute the difference Galois group
of the eigenvalue problem.
\end{abstract}
\section{Introduction}
We study the eigenvalue problem
$L\psi=\epsilon\psi$ for the ``$q$-Lam\'e'' difference operator $L$
in one complex variable $\lambda$:
\be
L\psi(\lambda)= 
\frac{\theta(\lambda-\gamma m,\tau)}
{\theta(\lambda,\tau)}
{\psi(\lambda+\gamma)}
+
\frac{\theta(\lambda+\gamma m,\tau)}
{\theta(\lambda,\tau)}
{\psi(\lambda-\gamma)}.
\ee
Here $\gamma$ is a complex parameter, $\tau$ is
a parameter in the upper half plane, and $m$ is
a nonnegative integer.

The theta function is Jacobi's $\theta_1$ function
\be
\theta(z,\tau)=-\!\sum_{j\in\Z+\frac12}e^{\pi ij^2\tau+2\pi ij(z+\frac12)}.
\ee
The difference operator $L$ appears in Sklyanin's work 
\cite{S} on
the representation theory of his elliptic algebra. In fact,
if we replace the theta functions in the numerators by the
other three Jacobi functions $\theta_\alpha$, we get difference
operators that obey with $L$ the commutation relations of 
Sklyanin's algebra.

The operator $L$ also appears as the simplest non trivial
quantum relativistic Ca\-lo\-ge\-ro--Moser Hamiltonian, introduced
by Ruijsenaars \cite{R}:
$L$ is conjugated to the
Ruijsenaars two-body Hamiltonian, and its 
semiclassical (or nonrelativistic) version is
the Lam\'e differential operator. More precisely,
as $\gamma\to 0$, $L=2{\mathrm{Id}}+{\mathrm{const}}\,
\gamma^2\ell
+O(\gamma^4)$, and $\ell$ is conjugated to
\be
-\frac {d^2}{d\lambda^2}+m(m+1)\wp(\lambda)+
{\mathrm{const}},
\ee
the Lam\'e differential operator of the classical
theory of ellipsoidal harmonics, see \cite{WW}. Here $\wp$ denotes
Weierstrass' doubly periodic function with periods
$1$ and $\tau$.

The eigenvalue problem for the Lam\'e operator 
has been studied in the classical literature. In
particular, Hermite gave a formula for eigenfunctions
which we would call today of ``Bethe ansatz'' type \cite{WW}:
he wrote a simple function depending on $m$
parameters, and proved that this function is
an eigenfunction if the parameters are solutions
of a set of $m-1$ algebraic equations (or $m$ equations if
one wants to specify the eigenvalue). He also
showed that all eigenfunctions with generic
eigenvalue are linear combinations of eigenfunctions
obtained  by this construction.

In more modern terminology, the Lam\'e operator
is an example of an ``algebraically integrable'' or
``finite gap'' Schr\"odinger
operator. The meaning of this is the content of
the following theorem, which follows from the
classical results, but whose true paternity
appears difficult to establish, see \cite{DMN}, \cite{CV}.

\begin{thm}\label{tuno}
There exists a differential operator $\ell_{2m+1}$ of order
$2m+1$, such that the algebra
$A$ of differential operators with meromorphic coefficients commuting
with $\ell$ is generated by $\ell$ and $\ell_{2m+1}$. Moreover $A$ is
commutative and $\ell\mapsto x$, $\ell_{2m+1}
\mapsto y$ defines an isomorphism from
$A$ to $\C[x,y]/(y^2-p(x))\C[x,y]$, for some polynomial $p$ of degree
$2m+1$.
\end{thm}

In other words the algebra of differential operators commuting with the
Lam\'e operator is isomorphic to the algebra of polynomial
functions on an affine hyperelliptic curve. This curve is called the
spectral curve of the differential operator $\ell$. Hermite's results
can be rephrased as saying that the spectral curve is
birational to a subvariety of the $m^{{\mathrm{th}}}$
symmetric power of the elliptic curve
$E_\tau=\C/\Z+\tau\Z$.

In \cite{KZ}, Krichever and Zabrodin considered the operator $L$,
and showed that eigenfunctions are parametrized by points on
a hyperelliptic curve. 

In \cite{FV}, we showed that $L$ is proportional to the transfer
matrix of the $2m+1$ irreducible representation of the elliptic quantum
group $E_{\tau,\gamma/2}(sl_2)$.  
The Bethe ansatz for this transfer matrix 
generalizes Hermite's results to the difference
case, see \cite{FV}, Theorem \ref{tFV}, and gives a formula for eigenfunctions
of $L$ parametrized by a Hermite--Bethe curve $Y$ given by $m-1$ equations
for $m$ variables on the elliptic curve $E_\tau$. In the formulation
of \cite{KZ}, we characterize double Bloch functions in terms
of their zeros. The results in \cite{KZ} and \cite{FV} are complementary:
an ``additive'' formula for eigenfunctions is given in \cite{KZ},
while a ``multiplicative'' formula is given in \cite{FV}. In both
cases the formula depend on parameters on a hyperelliptic curve.

In this paper, we study the difference operators
commuting with $L$. Our result is the following. Let $D_\gamma$ be the
algebra of difference operators of the form
\be
M\psi(\lambda)=\sum_{j=a}^{b} B_j(\lambda)\psi(\lambda+j\gamma),
\qquad a\leq b\in\Z
\ee
with 1-periodic meromorphic coefficients $B_j$. If $B_a$ and 
$B_b$ are non zero, we say that $M$ has {\it degree} $b$ and
{\it length} $b-a$. The degree and the length of a product is
the sum of the degrees and of the lengths of the 
factors, respectively.

\begin{thm}\label{tdue} Suppose that $\gamma$ is an irrational
number. Then 
there exists a difference operator $N$ of degree $2m+1$ and
length $4m+2$, such that the algebra $A$ of operators in $D_\gamma$
commuting with $L$ is generated by $L$ and $N$. Moreover $A$ is
commutative and $L\mapsto X$, $N\mapsto Y$ defines an
isomorphism from $A$ to $\C[X,Y]/(Y^2-P(X^2))\C[X,Y]$, for some
polynomial $P$ of degree $2m+1$.
\end{thm}
 
In the difference case, the spectral curve has two involutive automorphisms.
The automorphism $(X,Y)\mapsto (X,-Y)$ corresponds to the symmetry
$S\psi(\lambda)=\psi(-\lambda)$ of the eigenvalue problem, as in
the differential case. The automorphism $(X,Y)\mapsto (-X,-Y)$ 
corresponds to the involution $US\psi(\lambda)=
e^{\pi i\lambda/\gamma}\psi(-\lambda)$ which maps eigenfunctions
to eigenfunctions with opposite eigenvalue.

The operator $N$ can be written explicitly in terms of
 a remarkable family of difference operators $M_l$ commuting with
$L$:

\begin{thm}\label{ttre} Suppose that $\gamma\neq 0 \mod\Z+\tau\Z$.
For each generic $l\in\C$ there exists a difference
operator $M_l$ of the form
\be
M_l\psi(\lambda)=\sum_{k=0}^mA^l_{l-m+2k}(\lambda/\gamma)\psi(\lambda
+(l-m+2k)\gamma).
\ee
These operators commute pairwise and obey the relations
\be
LM_l=\frac {\theta(\gamma(l-m),\tau)}{\theta(\gamma l,\tau)}M_{l+1}
    +\frac {\theta(\gamma(l+m),\tau)}{\theta(\gamma l,\tau)}M_{l-1},
\ee
and the product rules
\be
M_lM_k=\sum_jA^l_j(k)M_{k+j}.
\ee
Moreover, $L$ is proportional to $M_{m-1}$ and $N$
is proportional to $M_{m+1}-M_{-m-1}$.
\end{thm}

We can then say, up to change of variable, that $M_l$ is an operator 
eigenfunction of the Lam\'e operator with ``eigenvalue'' $L$.

The eigenvalue problem $L\psi=\epsilon\psi$ becomes,
after the change of independent 
variable $\psi(\lambda)=u(\lambda)
\prod_{j=1}^m\theta(\lambda-j\gamma)
$,
\be
u(\lambda+\gamma)+
\frac
{\theta(\lambda+\gamma m,\tau)\theta(\lambda-\gamma(m+1),\tau)}
{\theta(\lambda,\tau)\theta(\lambda-\gamma,\tau)}u(\lambda-\gamma)
=\epsilon\, u(\lambda).
\ee
This equation has {\em elliptic} coefficients. Note incidentally
that if $\gamma\in\Z+\tau\Z$, the equation has constant coefficients
and can be easily solved in terms of exponential functions.

In the last section of this paper we compute the {\em difference
Galois group}, an analogue of the differential Galois
group of differential equations,
 of this equation in the case of irrational $\gamma$. It is
the Galois group of the extension of the difference
field of elliptic functions generated by the solutions.
It turns out that this Galois group is Abelian, in agreement
with the philosophy of \cite{BEG}, who characterized algebraically
integrable differential operators by the Abelian nature of
the corresponding differential Galois groups.

\section{Elliptic number notation}
We assume that $\gamma\neq 0\mod\Z+\tau\Z$.
It is convenient to introduce a new variable $x=\lambda/\gamma$. In
this variable, the step of the difference operator $L$ is integer,
and the periods are $\omega=\gamma^{-1}$ and $\omega'=\tau\gamma^{-1}$.
The ``elliptic number''
\be
[x]=\frac{\theta(\gamma x,\tau)}{\theta(\gamma,\tau)}
\ee
is an odd entire function of $x$ with zeros on the lattice
$\Z\omega+\Z\omega'$ and has transformation properties
\be
[x+\omega]=-[x],\qquad [x+\omega']=-e^{-\pi i\tau-2\pi i\gamma x}[x]
=-e^{-\frac{\pi i}{\omega}(\omega'+2x)}[x].
\ee
In this variable, which we use in the rest of this paper,
the $q$-Lam\'e operator reads
\be
L=\frac{[x-m]}{[x]}T_1+\frac{[x+m]}{[x]}T_{-1},\qquad T_j\psi(x)=\psi(x+j).
\ee
It will also be convenient to introduce ``elliptic binomial coefficients''
and ``elliptic factorials'':
\be
\left[{x\atop n}\right]=
\frac{[x][x-1]\cdots[x-n+1]}
{[n][n-1]\cdots[1]},\qquad 
\left[{x\atop 0}\right]=1,
\qquad [n]!=[1][2]\cdots [n].
\ee

\section{Bethe eigenfunctions}
We review the results of \cite{FV} on the $q$-Lam\'e operator. 
We assume that $\gamma\neq 0\mod\Z+\tau\Z$.
Let us introduce $m$ functions
of $t=(t_1,\dots,t_m)\in \C^m$:
\be
b_i(t)=\frac
{[t_i-m]}
{[t_i+m]}
\prod_{j:j\neq i}
\frac
{[t_j-t_i-1]}
{[t_j-t_i+1]}\,, \qquad i=1,\dots, m.
\ee
\begin{thm}\label{tFV}\cite{FV}
Let $(t_1,\dots,t_m,c)$ be a solution of the Bethe ansatz
equations:
\begin{equation}\label{baelame}
b_i(t)=
e^{2\gamma c},\qquad i=1,\dots,m,
\end{equation}
such that $t_i\neq t_j\mod \omega\Z+\omega'\Z$, if $i\neq j$. Then
\begin{equation}\label{psi3}
\psi(x)=e^{c\gamma x}\prod_{j=1}^m
{[x+t_j]},
\end{equation}
is a solution of the $q$-Lam\'e equation $\frac{[x-m]}{[x]}\psi(x+1)
+\frac{[x+m]}{[x]}\psi(x-1)=\epsilon_L\psi(x)$,
with eigenvalue
\be
\epsilon_L(t)=e^{-\gamma c}\frac{[2m]}{[m]}
\prod_{j=1}^m
\frac{[t_j+m-1]}
{[t_j+m]}.
\ee
\end{thm}
\noindent{\it Remark.\/} What was called $t_j$ in \cite{FV} is here
$(t_j-1/2)\gamma$. Also, $\gamma=2\eta$ in the notation of \cite{FV}.

Note that the functions $b_j$ have the propery that, for all $i,j$,
\bea
b_j(t_1,\dots,t_i+\omega,\dots,t_m)&=&b_j(t_1,\dots,t_i,\dots,t_m),\\
b_j(t_1,\dots,t_i+\omega',\dots,t_m)&=&
e^{4\pi i\gamma}b_j(t_1,\dots,t_i,\dots,t_m).
\eea
Thus if $(t_1,\dots,t_m,c)$ is a solution then  $(t_1,\dots
,t_i+\omega,\dots,c)$ and $(t_1,\dots,t_i+\omega',\dots,c+2\pi i)$
are also solutions. Moreover, the eigenfunctions are the same
for these solutions. Also, the equations and the eigenfunctions
are symmetric under permutations of the $t_j's$.  We thus have 
an action of the semidirect product $\Gamma=(\omega\Z+\omega'\Z)^m
\tilde\times S_m$ on $\C^m\times\C$, and the
Bethe ansatz equations take place in the quotient: The set
of solutions is
\bea
X_\infty&=&\{(t,c)\in (\C^m-D)\times\C
\,|\,b_j(t)=e^{2\gamma c}, j=1,\dots,m\}/\Gamma\,,\\
D&=&\cup_{i<j}\{t_i=t_j \mod \omega\Z+\omega'\Z\}.
\eea
On $X_\infty$ we have a $\Z$-action $(t,c)\mapsto (t,c+\pi i/\gamma)$.
The quotient is the set
\be
X=X_\infty/\Z
\ee
determined by the equations $b_i(t)=b_j(t)$. It is
 an algebraic subvariety of the symmetric power of 
our elliptic curve. The eigenvalue $\epsilon_L$ is a two-valued
function on $X$. It is a single-valued meromorphic function 
on the two-fold covering
\be
X_2=X_\infty/2\Z.
\ee
The eigenfunctions associated to points in $X_\infty$
which have the same projection on $X_2$ are proportional,
in the sense that they are related by multiplication
by a 
1-periodic function. 

In other words, a point in $X_2$ parametrizes a one-dimensional
vector space of eigenfunctions over the field $K$ of 1-periodic
meromorphic functions. Note that these eigenfunctions have all
the same eigenvalue.

Eigenfunctions associated to points in $X_2$ 
with the same projection on $X$
are related by the map $U:\psi\mapsto e^{\pi ix}\psi$. They have
opposite eigenvalues.

These varieties can be described rather explicitly in the neighborhood
of $c=\infty$:
\begin{lemma}\cite{FV}\label{lFV}
Let $\bar X$ be the closure of $X$ in the symmetric power $S^mE$ of
the elliptic curve $\C/\omega\Z+\omega'\Z$. Then $\bar X$ contains
the two  points $P_+=(-m+1,\dots,-1,0)$
and $P_-=(m-1,\dots,1,0)$. The irreducible component(s) of $X$ containing
 $P_+$ and $P_-$ form a curve $Y$ which is smooth at these two points.
In terms of a local coordinate vanishing at $P_+$, the eigenvalue
has the form
$\epsilon_L={\mathrm{const}}\,u^{-1/2}(1+O(u))$, and
$e^{2c\gamma}={\mathrm{const}}\,u^{-1}+O(1)$.
\end{lemma}

Let $Y_2\subset X_2$ be the double covering of $Y$. It is a curve
on which the eigenvalue is a single-valued function. The two
points $P_+$, $P_-$ are related by the symmetry $(t,c)\mapsto(-t,-c)$
of the Bethe ansatz equations.

\begin{thm}\label{tFV2}\cite{FV} Suppose that $\gamma\in\C$ is generic.
For generic $\epsilon\in\C$, there are precisely two  solutions
\be
(t_1,\dots,t_m,c)\qquad\mbox{and}\qquad
(-t_1,\dots,-t_m,-c), 
\ee
 of the Bethe ansatz
equations \Ref{baelame} in $X_\infty$ with given eigenvalue
$\epsilon_L=\epsilon$. The corresponding
eigenfunctions $\psi_\pm$ are linearly independent over the field $K$ of
$1$-periodic meromorphic functions of $x$, and all
solutions of the
$q$-Lam\'e equation $L\psi=\epsilon_L\psi$  are linear combinations of 
$\psi_+$, $\psi_-$ with coefficients in $K$. 
\end{thm}

\begin{definition}
The (many-valued) function on $X\times \C$ 
\be
\psi(t,x)=e^{c\gamma x}\prod_{j=1}^m\frac{[x+t_j]}{[t_j]}
\ee
where $c$ is determined by the Bethe ansatz equations,
is called the Baker--Akhiezer function. 
\end{definition}

The Baker--Akhiezer function is single-valued on $X_\infty\times\C$.

\begin{lemma}\label{lvan}
If $M=\sum_{j}A_j(x)T_j$ is a difference operator with
finitely many non-zero coefficients $A_j$ such that
$M\psi(t,\cdot)=0$ for all $t$ on the curve $Y$, then
$M=0$.
\end{lemma}
\begin{prf}
We consider the equation $M\psi(t,x)=0$ in the vicinity of the
point $P_+$ of Lemma \ref{lFV}: let $k$ be the largest number so that
$A_k\neq 0$. Since $\psi(t,x+j)/\psi(t,x)$ behaves as $e^{\gamma cj}
\sim u^{-j/2}$ as $u\to 0$, we see that for $M\psi/\psi$ to vanish,
it is necessary that its leading coefficient $A_k(x)$ vanishes,
a contradiction. Thus $M=0$.
\end{prf}

\section{Difference operators commuting with the $q$-Lam\'e
operator}\label{sfour}

We construct a sequence of difference operators that commute 
with the $q$-Lam\'e operator $L$. Let
\be
M_l=A^l_{l-m}(x)T_{l-m}+A^l_{l-m+2}(x)T_{l-m+2}+\cdots
+A^l_{l+m}(x)T_{l+m}
\ee
Usually $l$ is an integer, but we will occasionally take $l$ to
be a general complex number.
For $k=0,\dots, m$,
\bea
A^l_{l-m+2k}(x)
&=&
(-1)^k
\left[{m\atop k}\right]
  \prod_{j=0}^{m-k-1}    
   \frac{[l+m-j][x+m-j]}{[x+l+k-j]}\\
 & &\times
  \prod_{j=0}^{k-1}    
   \frac{[l-m+j][x-m+j]}{[x+l-m+k+j]},
\eea
with the understanding that a product over the empty set is one.

\begin{thm}
\
\begin{enumerate}
\item[(i)]
$M_{m}=\frac{[2m]!}{[m]!}{\mathrm{Id}}$, $M_{m-1}=\frac{[2m-1]!}{[m-1]!}L$.
\item[(ii)]
For all complex numbers $l$, $k$,
$M_lM_k=M_kM_l$.
\end{enumerate}
\end{thm}

\begin{prf}
(i) If $l=m$, all coefficient $A^l_{l-m+2k}$ with 
$k\neq 0$ vanish because of the factor $[l-m]$ in the second product.
If $k=0$,
\be
A^m_{0}(x)=
\left[{m\atop 0}\right]
  \prod_{j=0}^{m-1}    
   {[2m-j]}=\frac{[2m]!}{[m]!}
\ee
Similarly, if  $l=m-1$, the only non-vanishing coefficients
are $A^l_1(x)=[2m-1]![x-m]/[m-1]![x]=A^l_{-1}(-x)$.

\noindent (ii)
We first show that $M_l$ commutes with the $q$-Lam\'e operator
$L$. The equation $ML=LM$ for an operator of the form
$M=\sum A_jT_j$ is equivalent to the identities 
\be
A_j(x)\frac{[x\!+\!j\!+\!m]}{[x\!+\!j]}
\!+\!A_{j-2}(x)\frac{[x\!+\!j\!-\!2\!-\!m]}{[x\!+\!j\!-\!2]}
=
A_j(x\!-\!1)\frac{[x\!+\!m]}{[x]}
\!+\!A_{j-2}(x\!+\!1)\frac{[x\!-\!m]}{[x]},
\ee
for its coefficients. We have to show that our $A^l_j$, which are
zero except for $j$ in the set $\{l-m,l-m+2,\dots,l+m\}$,
are a solution of this equation. If we insert the formulae,
and cancel common factors, we see that the identity we have
to prove is
\bea
 &\displaystyle
\frac{[l\!-\!m\!+\!k\!-\!1][x\!-\!m\!+\!k\!-\!1][x\!+\!l\!+\!2k]}
{[k][x\!+\!l\!+\!k]}
-\frac{[l\!+\!k][x\!+\!k][x\!+\!l\!+\!2k\!-\!2m\!-\!2]}
{[m\!-\!k\!+\!1][x\!+\!l\!+\!k\!-\!m\!-\!1]}&
\\
 &=\displaystyle
\frac{[l\!-\!m\!+\!k\!-\!1][x\!+\!k]
[x\!-\!m\!-\!1][x\!+\!l\!+\!2k\!-\!m\!-\!1]}
{[k][x][x\!+\!l\!+\!k\!-\!m\!-\!1]}& 
\\
 &\displaystyle
-\frac{[l\!+\!k][x\!+\!m\!+\!1]
[x\!+\!k\!-\!m\!-\!1][x\!+\!l\!+\!2k\!-\!m\!-\!1]}
{[m\!-\!k\!+\!1][x][x\!+\!l\!+\!k]}.&
\eea
We make use of the properties of elliptic
numbers under translations by $\omega=\gamma^{-1}$ and
$\omega'=\gamma^{-1}\tau$.
Dividing both sides of this equation by $[x+k]$ yields
an equation, such that all terms are periodic functions
of $x$ with period $\omega$ and are multiplied by
$\exp(2\pi i\gamma(m-k+1))$, if $x$ is replaced by $x+\omega'$.
The (simple) poles are at $x=0,-l-k,-k$ and $-l+m-k+1$. However,
the difference between the left-hand side and right-hand side
has vanishing residue at these poles, as is easily checked,
for any generic value of $k$. But if $\alpha$ is a generic
complex number, the only  entire holomorphic function
$f$ such that $f(z+\omega)=f(z)$ and $f(z+\omega')=\alpha f(z)$
is $f=0$. Thus the identity is proved for generic, and,
by analiticity for all, values of $k$.  

To complete the proof of the commutativity of $M_l$ and $L$
we have to check the identity for the coefficients
 separately in the two extreme cases
$j=l-m$ and $j=l+m+2$, for which only two terms are non-zero.
This is easily done.

This completes the proof of (ii) if $l=m-1$. To prove the
general case, note that for any point $t\in X$,
\be
LM_l\psi(t,\cdot)=M_lL\psi(t,\cdot)=\epsilon_L(t)M_l\psi,(t,\cdot).
\ee
Hence $M_l\psi(t,\cdot)$ is an eigenfunction of $L$ with
eigenvalue $\epsilon_L(t)$. By Theorem \ref{tFV2}, it is proportional
to $\psi(t,\cdot)$. The same holds for $M_k$.
It follows that for generic $t$,
$[M_l,M_k]\psi(t,\cdot)=0$. Lemma \ref{lvan} implies then that
$[M_l,M_k]=0$.
\end{prf}

\begin{lemma}\label{l0}
Let $S$ and $U$ be the operators
$S\psi(x)=\psi(-x)$ and $U\psi(x)=e^{\pi ix}\psi(x)$.
Then, for all $l$, $M_lS=SM_{-l}$  and $UM_l=e^{\pi i(l-m)}M_lU$.
\end{lemma}
\begin{prf}
The first statement follows from the relation $T_jS=ST_{-j}$ and
the identity
\be
A^l_j(x)=A^{-l}_{-j}(-x).
\ee
The second follows from the relation $UT_{l-m+2k}=e^{\pi i(l-m)}T_{l-m+2k}U$.
\end{prf}

\begin{lemma}\label{l1}
Suppose that $\gamma$ is irrational and let $\omega=1/\gamma$.
Let $\phi(x)=\prod_{k=1}^m[x-k]$.
 Then, for any $j\in\C$, the first coefficient of any difference operator
$M=\sum_{k=0}^{l} B_{j-k}(x)T_{j-k}$ with $\omega$-periodic
coefficients $B_{j-k}$ which commutes with $L$ has
the form  $B_j(x)=c\frac{\phi(x)}{\phi(x+j)}$ for some constant $c$.
\end{lemma}
\begin{prf}
The $q$-Lam\'e operator $L$ can be written in the form
\be
L=\phi(x) T_1\phi(x)^{-1}+\phi(-x) T_{-1}\phi(-x)^{-1}.
\ee
It follows that $M_\phi=\phi^{-1}M\phi$ commutes with $T_1+C(x)T_{-1}$
for some $C$.
Thus the coefficient of $T_j$ in $M_\phi$ obeys the difference
equation $T_1f(x)=f(x)$ or $f(x+1)=f(x)$. But $f$ is also $\gamma^{-1}$%
-periodic since $M$ is in $D_\gamma$. If $\gamma$ is irrational,
it follows that $f$ is a constant $c$. Thus
$M=c\phi(x)T_j\phi(x)^{-1}+\cdots=c\phi(x)\phi(x+j)^{-1}T_j+\cdots$.
\end{prf}

\begin{lemma}\label{l4.5}
Let $\gamma$ be irrational. If $l\geq m+1$ or if $l\leq -m-1$ then
the  $M_l$ has  degree $l+m$ and length $2m$.
 If $l\in\{-m,-m+1,\dots,m\}$ then $M_l$ has degree $|m-l|$
and length $2|m-l|$.
\end{lemma}
\begin{prf}
This amounts to check when the coefficients $A^l_{k}$ vanish,
which is easy, given their factorized form.
\end{prf}

\begin{lemma}\label{l2}
Suppose that $\gamma$ is irrational and let $\omega=1/\gamma$.
Suppose that $j$ is either a generic complex number or a negative
integer. Then, any difference operator commuting with $L$
of the form
\be
M=B_j(x)T_j+B_{j-2}(x)T_{j-2}+\cdots+B_{j-2l}(x)T_{j-2l},\qquad B_j\neq 0,
\ee
with $\omega$-periodic coefficients $B_l$, has length $\geq 2m$.
\end{lemma}
\begin{prf} Let us assume that $M$ is nonzero and has length $<2m$.
The assumptions on $j$ ensure that $A^{j-m-2k}_{j-2k}(x)$ does not
vanish identically for $k=0,1,2,\dots$ 
Therefore we may subtract from $M$ a suitable linear combination
\be
M_c=c_0M_{j-m}+c_1M_{j-m-2}+\cdots+c_{m-1}M_{j-m-2(m-1)},
\ee
so that the result has degree $\leq j-2m$. We claim that 
$M=M_c$. Since $M_c$ has length $\geq 2m$ unless all $c_l$ vanish
(see Lemma \ref{l4.5}),
it then follows that $M=0$, contradiction. 
To prove our claim, let us suppose that
$M-M_c$  has degree $d\leq j-2m$. The
coefficient of $T_d$ in $M-M_c$, must have the form $a\,\phi(x)/\phi(x+d)$,
for some constant $a\neq 0$. In particular, it has a pole at $x=m-d$.
On the other hand, this coefficient is equal to
\begin{equation}\label{eqcA}
-c_0A^{j-m}_d(x)-\cdots-c_{m-1}A^{j-3m+2}_d(x),
\end{equation}
with no contribution from $M$, since $M$ has  length $<2m$. The
terms of this sum are of the form const $A^l_d(x)$ with $d\leq l+m-2$,
i.e., $A^l_{l-m+2k}(x)$ with $k\leq m-1$. The pole with largest
real part of $A^l_{l-m+2k}$ is $x=-l+m-k$, and this real part is
smaller than the real part of $m-d=-l+2m-2k$, if $k\leq m-1$. Therefore
\Ref{eqcA} is regular at $x=m-d$, a contradiction. 
\end{prf}

\begin{definition} Let $S$ be the involution $S\psi(x)=\psi(-x)$.
A difference operator $M$ is {\it symmetric} if $MS=SM$.
It is {\it antisymmetric} if $MS=-SM$.
\end{definition}

\begin{lemma}\label{l4}
Suppose that $\gamma$ is irrational.
Then any symmetric operator in $D_\gamma$ that commutes with $L$
is a polynomial in $L$.
\end{lemma}

\begin{prf}
Suppose that $M\in D_\gamma$ is a symmetric difference
operator of degree $j$ that commutes with $L$. Let us proceed
by induction. If $j=0$,
$M$ is a constant multiple of the identity by Lemma \ref{l1}.
Let $j>1$.
The coefficient of $T_j$ has then the form given in Lemma \ref{l1}.
If we subtract $c_jL^j$ from $M$ we obtain a symmetric difference operator of
degree at most $j-1$, which is a polynomial in $L$ by the induction
hypothesis.
\end{prf}

\begin{thm}\label{tquattro}
If $\gamma$ is irrational, then
all operators in $D_\gamma$ commuting with $L$ are polynomials
in 
$L$ and 
\be
N=M_{m+1}-SM_{m+1}S.
\ee
\end{thm}

\begin{prf}
We proved this for symmetric operators in Lemma \ref{l4}.
Assume that $M\in D_\gamma$ is antisymmetric and commutes with $L$.
If $M$ has degree $d\geq 2m+1$, we may subtract a constant multiple
of $NL^{d-2m-1}$, which is antisymmetric and
has degree $d$ by Lemma \ref{l4.5}, to get 
an operator of smaller degree. This is possible by Lemma \ref{l1}.
Thus we may assume that $M$ is antisymmetric and 
has degree $\leq 2m$. By subtracting from $M$ a suitable 
polynomial of $L$ of degree $\leq 2m$, we get an operator
commuting with $L$ and
of the form $M'=B_{-1}(x)T_{-1}+\cdots+B_{-2m}(x)T_{-2m}$.
To apply Lemma \ref{l2}, we write $M'=M'_o+M'_e$, where
$M'_o=B_{-1}T_{-1}+B_{-3}T_{-3}+\cdots$ is the sum of the odd terms.
Both $M'_o$, $M'_e$ commute then with $L$, since the property
of commutation with $L$ is equivalent to relations involving
only even or odd coefficients.
But both $M'_o$ and $M'_e$ are of negative degree and length 
$< 2m$, and therefore vanish\
by Lemma \ref{l2}.
\end{prf}

Thus the algebra $A$ is generated by $L$ and $N$. Since
$N$ is antisymmetric, its square is symmetric and by
Lemma \ref{l4} we have a relation
\be
N^2=Q(L)
\ee
for some polynomial $Q$. Comparing the degrees we see that $Q$
has degree $4m+2$. By Lemma \ref{l0}, $N^2$ commutes with
$U$, but $LU=-UL$. Since the powers $L^j$ of $L$ are linearly 
independent (they have different degree), it follows that 
$Q(L)=P(L^2)$ for some polynomial $P$ of degree $2m+1$.

We are now ready to complete the proof of Theorem \ref{tdue}.
The fact that $L$ and $N$ commute and obey this relation
means that $X\mapsto L$, $Y\mapsto N$ defines
a surjective homomorphism of algebras
\be
h:\C[X,Y]/(Y^2-P(X^2))\C[X,Y]\to A
\ee
Any element of the left algebra is represented uniquely by
a polynomial of the form $f(X)+g(X)Y$. Such a polynomial is
in the kernel of $h$ if and only if $f(L)+g(L)N$ vanishes. But
this means that the symmetric and antisymmetric
parts $f(L)$ and $g(L)N$ vanish separately. By considering
the coefficient of $T_j$ with highest $j$, we see inductively that
all coefficients of the polynomial $f$ and $g$ vanish. Thus
$\phi$ is an isomorphism.

\section{Eigenvalues}
In this section we compute the eigenvalues of our commuting operators
on the Baker--Akhiezer function  $\psi(t,x)$. This eigenvalue map
maps an element of the algebra of commuting difference operators
(i.e., a function on the hyperelliptic curve)
to a two-valued function on the 
Hermite--Bethe curve $Y$ and, as will be shown,
realizes the birational equivalence between the hyperelliptic
curve and a double covering of the Hermite--Bethe curve.

We start by describing some remarkable properties of the
difference operators $M_l$, $l\in\C$.

\begin{proposition}
For all generic complex $l$,
\be
LM_l=
\frac{[l+m]}{[l]}M_{l-1}+
\frac{[l-m]}{[l]}M_{l+1}
\ee
\end{proposition}
\begin{prf}
If $l$ is generic, then
the operator $LM_l$ is of degree $l+m+1$ and length $2m+1$ (Lemma
\ref{l4.5}). Therefore,
by Lemma \ref{l1}, 
\be
LM_l=C_l\frac{\phi(x)}{\phi(x+l+m+1)}T_{l+m+1}+\cdots,
\ee
up to terms of lower degree. Here $C_l$ appears in the coefficient
of $T_{l+m}$ in $M_l$: $A^l_{l+m}(x)=C_l\phi(x)/\phi(x+l+m)$.
It follows that by subtracting a suitable multiple of $M_{l+1}$
from $LM_l$ we get an operator of degree $\leq l+m-1$. Similarly,
$LM_l=SLM_{-l}S=C_{-l}\phi(-x)/\phi(-x-l+m+1)T_{l-m-1}$ plus terms
of higher degree, and we may kill the coefficient of $T_{l-m-1}$
by subtracting a multiple of $M_{l-1}$. We conclude that
\be
LM_l-\frac{C_l}{C_{l+1}}M_{l+1}-\frac{C_{-l}}{C_{-l+1}}M_{l-1}
\ee
is an operator of length $<2m$ commuting with $L$, and thus 
vanishes by Lemma \ref{l2}. The ratios of $C_l$ can easily
be computed from the explicit expression for $A^l_{l+m}$,
and give the desired result.
\end{prf}

Thus $M_l$, viewed as a function of $l$, is an eigenvector of
the $q$-Lam\'e operator in the space of difference operators
with ``eigenvalue'' $L$.
\begin{lemma}\label{ll2}
Let $\omega=1/\gamma$. Then $M_{l+\omega}=M_lT_\omega$.
\end{lemma}
\begin{prf}
We have 
\bea
M_{l+\omega}&=&\sum_{k=0}^mA_{l+\omega-m+2k}^{l+\omega}
T_{l+\omega-m+2k}\\
 &=& 
\sum_{k=0}^mA_{l-m+2k}^{l}
T_{l-m+2k}T_\omega,
\eea
since $A_{l-m+2k}^{l}$ is $\omega$-periodic as a function of $l$.
\end{prf}

\begin{proposition}\label{p5.3}
\be
M_l\psi(t,x)=\epsilon_l(t)\psi(t,x),
\ee
with eigenvalue
\be
\epsilon_l(t)=\frac{[2m]!\psi(t,l)}{[m]!\psi(t,m)}\, ,
\ee
\end{proposition}
\begin{prf}
Let $\psi(t,x,l)=M_l\psi(t,x)$ and denote $\hat L$ the $q$-Lam\'e
operator acting on the variable $l$. We then have
\bea
\hat L\psi(t,x,l)
&=&
\frac{[l+m]}{[l]}M_{l-1}\psi(t,x) 
+\frac{[l-m]}{[l]}M_{l+1}\psi(t,x) \\
&=& LM_l\psi(t,x) \\
&=& M_lL\psi(t,x) \\
&=& \epsilon_L(t)M_l\psi(t,x).
\eea
Let $\omega=1/\gamma$. Let $e^c$ be the multiplier of $\psi(t,x)$:
$\psi(t,x+\omega)=e^c\psi(t,x)$. Then, by Lemma \ref{ll2},
we have $\psi(t,x,l+\omega)=e^c\psi(t,x,l)$. In other words, both
as function of $x$ and as a funcion of $l$,
$\psi(t,x,l)$ is an eigenfunction of $L$ with the same eigenvalue
and multiplier. If $t$ is generic, there is only one such eigenfunction
up to normalization. Thus
\be
\psi(t,x,l)=f(t)\psi(t,l)\psi(t,x),
\ee
for some $f(t)$. On the other hand, we know that if $l=m$,
$M_m=[2m]!/[m]!$ times the identity. This determines $f$ and
we get
\be
\psi(t,x,l)=\frac{[2m]!}{[m]!}\frac{\psi(t,l)\psi(t,x)}
{\psi(t,m)}\,.
\ee
\end{prf}

As a corollary, we see that the relation of Theorem \ref{ttre}
is a special case of more general product rules:

\begin{corollary}\label{cmlmk}
For generic $l,m\in\C$,
\be
M_lM_k=\sum_jA^l_j(k)M_{k+j},
\ee
where $A^l_j$, $j=l-m,l-m+2,\dots,l+m$,
 are the coefficients defined in Section \ref{sfour}.
\end{corollary}
\begin{prf} By Lemma \ref{lvan}, it is sufficient to prove this
identity for the eigenvalues $\epsilon_l(t)$. We first note that
Proposition \ref{p5.3} can be rewritten, after replacing $x$ by $k$, as
\be
\sum_jA^l_j(k)\psi(t,k+j)=\epsilon_l(t)\psi(t,k).
\ee
Hence,
\bea
\sum_jA^l_j(k)\epsilon_{k+j}(t)&=&
\frac{[2m]!}{[m]!}\sum_jA^l_j(k)\frac{\psi(t,k+j)}{\psi(t,m)}
\\
 &=&\frac{[2m]!}{[m]!}\epsilon_l(t)\frac{\psi(t,k)}{\psi(t,m)}
\\
 &=&\epsilon_l(t)\epsilon_k(t).
\eea
\end{prf}

\begin{corollary}
Let $N=M_{m+1}-M_{-m-1}$. Then
\be
N\psi(t,x)=\epsilon_N(t)\psi(t,x),
\ee
where
\be
\epsilon_N(t)=\frac{[2m]!}{[m]!}
\biggl(e^{\gamma c}\prod_{j=1}^m\frac{[m+t_j+1]}{[m+t_j]}
-e^{-\gamma c}\prod_{j=1}^m\frac{[m-t_j+1]}{[m-t_j]}\biggr),
\qquad e^{\gamma c}=\sqrt{b_j(t)}.
\ee
The map $t\mapsto(\epsilon_L(t),\epsilon_N(t))$ defines a
birational isomorphism from the double covering $Y_2$ of the 
Hermite--Bethe curve to
the curve $\{(X,Y)\in\C^2| Y^2=P(X^2)\}$.
\end{corollary}
\begin{prf} The expression for the eigenvalue is taken from 
Proposition \ref{p5.3}.

By construction, the 
function $\epsilon_L(t)$ is a two-to-one rational function from
the closure of $Y_2$
onto $\Pee^1$. The two points in $\epsilon_L^{-1}(p)$ for 
generic $p\in\Pee^{1}$
are related by the symmetry $(t,c)\mapsto (-t,-c)$. The eigenvalue
$\epsilon_N$ is odd under this symmetry, and we thus have a
one-to-one  (at generic points) map from $Y_2$ to the hyperelliptic curve.
\end{prf}

{\it Remarks.\/} 
\begin{enumerate}
\item[1.] The hyperelliptic 
curve $Y^2=P(X^2)$ has a double point at infinity
which is resolved into the two points $P_+$ and $P_-$ of the double
covering $Y_2$ of the Bethe--Hermite curve. These are the points at
which the eigenvalues of the commuting operators have poles.  Our
results can be considered as
a degenerate case of the difference version of Krichever's
construction: to a smooth projective curve $C$ of genus $g$, two
points $P_+,P_-$ on it and a generic effective
divisor $D$ of degree $g$, Krichever \cite{K} associates a
Baker--Akhiezer function $\psi(p,x)$, a properly normalized
 meromorphic function of $p\in
C$ and $x\in\Z$ with divisor (the formal sum-with-multiplicities
of the poles minus the zeros)
 $x(P_+-P_-)+D$. To each meromorphic function
$f$ on $C$ which is regular on $C-\{P_+,P_-\}$ is associated a
difference operator $M_f$ with integer steps,
for which $\psi(p,\cdot)$ 
is an eigenfunction with eigenvalue $f(p)$, and $f\mapsto M_f$ is
an algebra homomorphism. In particular, if $P$ is a polynomial of
odd degree $2m+1$ without
multiple roots, then the curve $Y^2=P(X^2)$ can be compactified to
a smooth hyperelliptic curve by adding two points $P_+$, $P_-$ at infinity,
and one has, upon choosing a divisor, a pair of commuting 
difference operators $M_X$ and $M_Y$. Note that in general the
Baker--Akhiezer function is written in terms of the Riemann theta
function of the hyperelliptic curve, whereas in this case it
can be written purely in terms of an elliptic curve.

\item[2.] Note the analogy with ``fusion algebras''. Let
$N^i_{jk}=\dim {\mathrm{Hom}}(V_i,V_j\otimes V_k)$ be the dimensions
of the space of homomorphisms of, say, representations of a simple
Lie group $G$. 
Here $V_i$ are irreducible finite dimensional representations, labeled
by their highest weight $i\in P^+$. Let us introduce difference 
operators $M_l$ acting on functions on $P^+$ by the formula
$M_lf(i)=\sum_jN^i_{lj}f(j)$. The commutativity and
associativity of the tensor
product imply that the operators $M_l$ commute with each other.
Moreover, we have $M_lM_k=\sum_jN^j_{lk}M_j$, cf.~Corollary
\ref{cmlmk}. Let $\psi(t,j)=\tr_{V_j}(t)$
be the character of the representation $V_j$. It is a function of
$t\in T/W$, the quotient of a Cartan torus by the Weyl group. Then,
for fixed $t$, the function $\psi(t,j)$ of $j$ is a common
eigenfunction of all difference operators $M_l$:
\be
M_l\psi(t,\cdot)=\epsilon_l(t)\psi(t,\cdot),
\qquad \epsilon_l(t)=\frac{\psi(t,l)}
{\psi(t,0)},
\ee
cf.~Proposition \ref{p5.3}.
If we replace $G$ by a quantum group at root of unity, the same
formulae apply, except that $P^+$ is replaced by a finite subset and
$T$ is replaced by the set of points of a finite order $N$, depending
on the order of the root of unity, in the Cartan torus. In this case,
$\psi(t,j)$ has a remarkable interpretation in terms of
representations of $SL(2,\Z)$, discovered by E. Verlinde \cite{V}, in the
context of conformal field theory.  Is there a similar interpretation
in our case?
\end{enumerate}

\section{The Galois group}
We compute the ``difference Galois group'' of the $q$-Lam\'e
equation. This group is a difference analogue of the
differential Galois group of differential equations. The computation
is motivated by the recent result of Braverman, Etingof and
Gaitsgory \cite{BEG}, who, in the differential case,
 characterized algebraic integrability
by the Abelian property of the Galois group.

\begin{definition}
A {\em difference field} is a field $F$ together with an automorphism
$T\in{\mathrm{Aut}}(F)$. An {\em extension} $E\subset F$
of difference fields is a subfield $E$ of a difference field $F$
such that $T(E)\subset E$. An {\em automorphism} of a difference
field $F$ is an automorphism of the field $F$ commuting with $T$.
The {\em Galois group} of an extension $E\subset F$ of difference
fields is the group of those automorphisms of $F$ which restrict to
the identity on $E$.
\end{definition}

Let us now consider our difference equation $L\psi=\epsilon\psi$
for fixed generic $\epsilon$. If we make the change of variables
$\psi(x)=u(x)\prod_{j=1}^m[x-j]$, the equation becomes
\begin{equation}\label{eu1}
u(x+1)+\frac{[x+m][x-m-1]}{[x][x-1]}u(x-1)=\epsilon u(x).
\end{equation}
This    difference equation has  coefficients  in the field
$E$ of elliptic functions with periods $\omega,\omega'$. 
The field $E$ is a difference field with $T$ the shift $Tf(x)=f(x+1)$.

Let $F$ be the differential field generated over $E$ by the
meromorphic solutions of the difference equation \Ref{eu1}.
It is the field of all rational functions in the solutions
and their images by $T^j$, $j\in\Z$,
with coefficients in $E$.

It follows from Theorem \ref{tFV} that the solutions of \Ref{eu1}
are linear combinations of $u_+(x)$ and $u_-(x)=u_+(-x)$ with
coefficients in $K$, the field of meromorphic $1$-periodic
functions. The solution 
$u_+$ has the form
\be
u_+(x)=e^{c\gamma x}\prod_{j=1}^m\frac{[x+t_j]}{[x-j]}\,.
\ee
\begin{thm}\label{tgal}
 Let $\gamma\in\R-\Q$ and $\epsilon$ be generic. Then the
Galois group of the extension $E\subset F$ of difference fields
is isomorphic to  the group $K^\times$ of non-zero meromorphic
1-periodic functions on the complex plane. A function $h\in K^\times$
corresponds to the automorphism acting on solutions by
$u_\pm(x)\mapsto h(x)^{\pm1}u_\pm(x)$.
\end{thm}

To prove the theorem, we first need some auxiliary results.

\begin{lemma}\label{lg1} If $\Phi$ is an element of the Galois
group then $\Phi(u_\pm)=h^{\pm 1}u_\pm$, for some $h\in K^\times$.
\end{lemma}

\begin{prf}
The function $u_+(x)$ has  constant multipliers as $x$ is
shifted by $\omega$ or $\omega'$. Thus $Tu_+/u_+$ is an elliptic
function.
It follows that if a Galois automorphism sends $u_+$ to a function
$\tilde u_+$, then $T\tilde u_+/\tilde u_+=Tu_+/u_+$. Thus
$\tilde u_+=hu_+$ for some function $h\in K^\times$.

Similarly, $u_+(x)u_-(x)$ is elliptic. Thus if $u_+$ is sent to
$h u_+$, then $u_-$ is sent to $h^{-1}u_-$.
\end{prf}
 Let $\hat E$ be the field generated by $K$
and $E$.
\begin{lemma}\label{lg2} $F=\hat E(u_+)$, i.e.,
the field $F$ consists of rational functions in $u_+$ with 
coefficients in $\hat E$.
\end{lemma}
\begin{prf} Notice first that $K\subset F$: an element $h\in K$
is the ratio of solutions $hu_+/u_+$. It remains to show that
every element of $F$ can be written as a rational function in $u_+$
with coefficients in $\hat E$. But the proof of the preceding lemma
shows that $u_-$ and $Tu_+$ and thus all solutions, as well as
all their images by $T^j$, $j\in\Z$, are rational
functions in $u_+$ with coefficients in $\hat E$.
\end{prf}
\begin{lemma}\label{lg3}
$u_+$ is transcendental over $\hat E$.
\end{lemma}
\begin{prf}
This means that  $u_+$ is not the solution
of any non-trivial polynomial equation with coefficients in
$\hat E$. Suppose that there is such an equation
\be
P(x,u_+(x))=\sum A_j(x)u_+(x)^j=0,\qquad A_j\in \hat E.
\ee
Since $\omega$ is real, there exists a stricly increasing
sequence of integers $n_1, n_2,\dots$, so that
the distance between $n_l\omega$  
and the lattice of integers converges to zero. Let 
$x$ be any generic complex number and set $x_l=x+n_l\omega$.
Then $\lim_{l\to\infty}A_j(x_l)=A_j(x)$, for all $j$.
On the other hand, $u_+(x_l)=C^{n_l}u_+(x)$, for some non-trivial
constant $C$. Since $P(x_l,u_+(x_l))=0$ for all $l$, it
follows that all coefficients $A_j(x)$ vanish at $x$. But
$x$ is arbitrary. Therefore all coefficients
vanish identically, a contradiction.
\end{prf}

The proof of Theorem \ref{tgal} can now be completed. 
What is left to prove is that, for every $h\in K^\times$,
there exists a unique Galois automorphism sending
$u_+$ to $hu_+$. 

The uniqueness follows from Lemma \ref{lg1} and
Lemma \ref{lg2}.

To prove existence, we have to show that for all rational
functions $f\in\hat E(X)$ of one indeterminate, the map
\begin{equation}\label{egalois}
f(u_+)\mapsto f(hu_+),
\end{equation}
is well-defined, i.e., independent of the choice of the
function $f$ used to represent an element $f(u_+)\in F$,
and that it defines an automorphism of difference fields.

The map is well-defined: if $f(u_+)=g(u_+)$ for rational
functions $f=p/q$, $g=r/s$, then $u_+$ is a solution
of the polynomial equation $ps-qr=0$. By Lemma \ref{lg3}, this
equation must be trivial, meaning that $f=g$ in $\hat E(X)$.
In particular $f(hu_+)=g(hu_+)$.

It is clear that \Ref{egalois} defines an automorphism 
of fields with inverse $f(u_+)\mapsto f(u_+/h)$. Let us
show that it is an automorphism of difference fields.
We have $Tf(u_+)=\bar f(u_+)$. The rational function
$\bar f(X)$ has the form $(Tf)(a X)$, where 
$a=Tu_+/u_+\in E$ and $Tf$ is obtained
from $f$ by acting with $T$ on the coefficients.
We have to show that $T(f(hu_+))=\bar f(hu_+)$. But
 since $T(hu_+)=hT(u_+)$, we have
\be
T(f(hu_+))=(Tf)(hT(u_+))=(Tf)(h\, a\, u_+)=\bar f(hu_+).
\ee
The proof is complete.

\noindent{\it Remark.\/} Our construction is a special case of
a more general construction in higher dimension: a difference
field in $n$ dimensions is a field $F$ together with $n$ commuting
automorphisms $T_1,\dots,T_n$. Extensions and Galois groups are
defined as obvious generalizations of the $n=1$ case.
Suppose that $F$ is a field of functions 
$f(x_1,\dots,x_n)$ and $T_i$ are shift operators 
\be
T_i f(x_1,\dots,x_n) = f(x_1,\dots,x_i+a_i,\dots,x_n),\qquad i=1,\dots,n.
\ee
Let us say that the function $g$ is {\em elementary}
 with respect to the difference field $F$
if there are functions $f_1,\dots,f_n$ in $F$ such that $g$ satisfies
the equations
\be
T_i g / g = f_i,\qquad    i=1,\dots,n.
\ee
Let  
\be
L_j F(x_1,\dots,x_n)=0,\qquad  j=1,\dots,k,
\ee
be a system of linear difference equations with coefficients in the difference
field $F$. Assume that the space of solutions has a basis consisting of
elementary functions. Then the Galois group of this system is Abelian.

\end{document}